\documentclass[12pt,reqno]{article}
\usepackage{amscd,  amsmath, amssymb, amsthm, sectsty, url, enumerate}

\newtheorem{thm}{Theorem}[section]

\newtheorem{lem}[thm]{Lemma}
\sectionfont{\normalsize}

\theoremstyle{remark}

\theoremstyle{definition}

\numberwithin{equation}{section}

\title{MIND SWITCHES IN FUTURAMA AND STARGATE}
 
\author{\\ Ron Evans\\
Department of Mathematics\\
University of California, San Diego\\
La Jolla, CA  92093-0112\\
revans@ucsd.edu
\\ \\
Lihua Huang\\
Department of Mathematics\\
Rutgers University\\ 
Piscataway, NJ 08854-8019\\
lhuang@math.rutgers.edu
\\ \\ 
} 
\date{August 2014}
\begin{document}
\maketitle
\thispagestyle{empty}
\newpage
\setcounter{page}{1}
\begin{abstract}

We generalize two mind-switching problems that arise in connection
with the popular sci-fi television series {\em Futurama} and
{\em Stargate SG-1}.
Optimal solutions to these problems are found by 
answering the following question about a permutation $\sigma$
expressed as a product of nontrivial disjoint cycles.
``When writing $\sigma$ as a product of distinct transpositions
none equal to a factor of $\sigma$, what is the smallest number 
of transpositions that can be used?"

\end{abstract}

\section{Introduction: Problems from Futurama and Stargate}

``The Prisoner of Benda" \cite{ams}, \cite{benda}, 
a 2010 episode of
{\em Futurama}, features a two-body mind-switching machine.
Any pair can enter the machine to swap minds, but the machine
has the limitation that it will not work more than once on the
same pair of bodies. A two-body mind-switching machine with
exactly the same limitation is featured in
``Holiday" \cite{holiday},
a 1999 episode of {\em Stargate SG-1}.
In both episodes, the participants in the mind swapping
want to return to their original bodies, but they don't know how.  Brilliant
characters (Sweet Clyde Dixon in {\em Futurama}, Captain Samantha Carter in
{\em Stargate SG-1}) save the day by figuring out how to reverse the switching.

We first present a problem arising from
{\em Futurama}. In ``The Prisoner of Benda,"
Amy and the Professor enter the machine and swap minds.
Unhappy with their new bodies,
they re-enter the machine and attempt to undo the switch.
When nothing changes, 
the Professor recalls that
due to cerebral immune response, ``Once two bodies have switched minds,
they can never switch back."  
He proposes the idea of employing a third body
for temporary mind storage, 
but he soon realizes that three bodies
will not suffice to restore all minds to their original bodies.
Agitated, Amy exclaims, ``Oh  no! Is it possible to get everyone
back to normal using four or more bodies?"  The Professor responds,
``I'm not sure.  I'm afraid we need to use...math!"

Throughout the episode, nine characters take
part in a mind-swapping spree: 
Fry, Zoidberg, Amy, Hermes,
Leela, Professor, Bender, Emperor, and Washbucket.
For brevity, call their bodies
$1,2,3,4,5,6,7,8$, and $9$, respectively. 
These characters engage 
in a sequence of seven swaps:
first the pair 3,6 swap minds, then 3,7, then 5,6, then 3,9, then 1,2, then
8,9, and finally 4,5.  

After many adventures in each other's bodies, the characters 
decide that they were better off before the switching.
The Professor (in Bender's body) tries in vain to 
find a way to restore all minds to their original bodies.  
Sweet Clyde Dixon comes to the rescue with his proof of an inversion theorem
displayed on a greenboard.  The theorem shows how to undo
the switching with a sequence of thirteen swaps, as given in
\eqref{q13} below.
Sweet Clyde becomes a hero and the Emperor makes him a duke.
Curiously though, the characters actually restore normalcy with a sequence
of thirteen swaps {\em different} from that shown on the greenboard; first
$s$ and $1$ swap minds, then the pair $t$,2, then $s$,2, then $t$,1,
then $s$,6, then $t$,9, then $s$,5, then $t$,8, then $s$,4, then $t$,7,
then $s$,3, then $t$,6, and finally $s$,9, where $s$ and $t$ denote
the bodies of Globetrotters Sweet Clyde and Bubblegum Tate.

In a popular 2010 video \cite{grime}, 
Cambridge University mathematician James Grime
points out that the switching could be reversed without employing
the bodies $s$ and $t$, using the following nine swaps:
first $2$ and $3$ swap minds, then the pair 1,9, then 1,8, then 1,7, then 1,6,
then 1,5, then 1,4, then 1,3, and finally 2,9.  
Grime proceeds to ask if 9 is the {\em smallest}
number of swaps that could accomplish the reversal.
In the next section, this question is reformulated using mathematical 
terminology.
Then in Theorem 3.4, we give an affirmative answer to a 
general version of the question
involving an arbitrary number of bodies.

We next turn to a problem arising from {\em Stargate SG-1}.
In ``Holiday,"  
a crisis is created when the aging Machello tricks young handsome Daniel
into swapping minds with him. As Jack and
Teal'c attempt to help their colleague Daniel
by retrieving the mind-switching machine,
they blunder into swapping themselves.
Unaware that Machello has designed the machine to block a direct reversal,
physicist Samantha Carter attempts in vain to rectify the blunder by 
instructing Jack and Teal'c to grasp the machine's handles in different ways.
Meanwhile, Machello is enjoying his holiday in Daniel's body.  

Machello later
feels guilty about what he has done to Daniel, and tells Samantha,
``If I could trade places with him again, I would.  
But unfortunately, I am the only one who cannot."
His last sentence sparks Samantha's idea on how to reverse the
switching.  She tells Jack, ``We're going to have to play a little
musical chairs with your bodies."  At Samantha's instruction, first
the bodies of Teal'c and Daniel swap minds, then Jack and Machello,
then Teal'c and Machello, and finally Jack and Daniel.
This sequence of four swaps brings everyone back to normal, and 
Samantha earns the undying gratitude of the Stargate crew.
Is 4 the {\em smallest} number of swaps that can accomplish the reversal?
In Theorem 3.5, we give an affirmative
answer to a general version of this question
involving an arbitrary even number of bodies.

\section{Mind switching and permutations}

In this section, we generalize the aforementioned problems and 
frame them in the language of group theory.  Let $S_n$ denote the
group of $n!$ permutations of $\{1,2,\dots, n\}$.

We begin with the problem from {\em Futurama}.
Recall that the nine characters engage in a sequence of seven swaps,
beginning with the pair 3,6.
We represent this sequence by the 
right-to-left formal product
\begin{equation}\label{p9}
B_9: = (45)\cdot(89)\cdot(12)\cdot(39)\cdot(56)\cdot(37)\cdot(36),
\end{equation}
where the factors in \eqref{p9} are
transpositions ($2$-cycles) in the symmetric group $S_9$.

In the sequel, the symbol $B$ (with or without a subscript) will 
always denote a formal product of distinct transposition factors in $S_n$,
representing a sequence of swaps that one aims to reverse.
For brevity,  ``transposition factors" will simply be called ``factors".
When the factors of $B$ are written in reverse order, the resulting
formal product will be denoted by $B^{-1}$.  The product $B$ is an
element in the group of words on the alphabet of transpositions.

One can view $B$ as a permutation in $S_n$ via
right-to-left multiplication in $S_n$ of its factors.
This permutation will be denoted by $\sigma(B)$.  
Every permutation in $S_n$ can be expressed as a product of
nontrivial disjoint cycles \cite[Theorem 2.3.5]{BB}, and we will
always assume that $\sigma(B)$ is expressed in this form.
For example, the product $B = (12)\cdot(23)$ effects the $3$-cycle
$\sigma(B)=(123) \in S_3$.  This 3-cycle represents the permutation 
that sends 1's mind to 2, 2's mind to 3, and 3's mind to 1.
From \eqref{p9},
\begin{equation}\label{dp9}
\sigma(B_9) = (12)(3456789).
\end{equation}

The mind-switching history is encoded in the product $B$, but not generally
in the permutation $\sigma(B)$.  To illustrate this, consider the two products
$B=(12)\cdot(23)$ and $B'=(23)\cdot(13)$.  
The product $B$ tells us that the pair 1,2
swapped right after the pair 2,3, while $B'$ tells us that 2,3 swapped
right after 1,3.
Although $B$ and $B'$ 
are completely different formal products, they effect the
same permutation $\sigma(B)=\sigma(B')=(123)$.
This permutation conveys no information that distinguishes 
between the distinct switching histories encoded by $B$ and $B'$.

The participants in {\em Futurama}'s 
mind-swapping spree wish to undo $B_9$; that is, they 
want to reverse the switching 
in order to bring all minds back to their original bodies.
To undo a general product of transpositions $B$
representing a sequence of mind swaps,
one needs to construct a formal product $Q$ of distinct transpositions such that
\begin{enumerate}[(i)]
\item the products $Q$ and $B$ have no factor in common; and
\item  $\sigma(Q)=\sigma(B)$.
\end{enumerate}
Then $\sigma(Q^{-1}\cdot B)$  equals the identity permutation and 
$Q^{-1}$ undoes $B$.
To reverse the switching as efficiently as possible,
we want to undo $B$ using a product $Q$ that has the 
smallest possible number of factors. Any such $Q$ is said to be ``optimal" or
``best possible".  In the sequel, the symbol $Q$ (with or without a subscript) 
will always denote a formal product of distinct transpositions.

For the product $B_9$ in \eqref{p9}, which represents the 
sequence of mind swaps in ``The Prisoner of Benda,"
Grime \cite{grime} observed that
$\sigma(B_9)=\sigma(Q_9)$, where
\begin{equation}\label{q9}
Q_9:=
(23)\cdot(19)\cdot(18)\cdot(17)\cdot(16)\cdot(15)\cdot(14)\cdot(13)\cdot(29).
\end{equation}
The product $Q_9$ has no factor in common with the product $B_9$.
Thus $Q_9^{-1}$ serves to undo $B_9$ using nine transpositions.
Grime asked if $Q_9$ is best possible.  We will give an
affirmative answer to a more general question.

For $n \ge 5$, let $B_n$ be any product of distinct transpositions
representing a sequence of mind swaps with the following two properties:
\begin{enumerate}[(i)]
\item bodies 1 and 2 swap with each other but not
with any other body,  as in \eqref{p9}; and 
\item  as in \eqref{dp9},
\begin{equation}\label{dpn}
\sigma(B_n) = (12)(345\dots n).
\end{equation}
\end{enumerate}
Then analogous to the case $n=9$ discussed in the preceding paragraph, we have
$\sigma(B_n)=\sigma(Q_n)$, where
\begin{equation}\label{qn}
Q_n:=
(23)\cdot\left((1n)\cdots(15)\cdot(14)\cdot(13)\right)\cdot(2n).
\end{equation}
The product $Q_n$ has no factor in common with the product $B_n$.
Thus $Q_n^{-1}$ serves to undo $B_n$ using $n$ transpositions.
Is $Q_n$ best possible? 
For $n=9$, this is Grime's question.
Theorem 3.4 gives an affirmative answer for each $n \ge 5$.

``The Prisoner of Benda" was written by Ken Keeler \cite[p. 207]{singh}, 
who earned
a Ph.D. in applied mathematics at Harvard University in 1990.
For the show, Keeler presented an algorithm  \cite{benda}
that is designed to 
reverse any sequence of swaps, i.e., it is designed to undo any product $B$.
We gave a best possible algorithm 
optimizing Keeler's in \cite[Theorem 1]{keeler}.
Both algorithms depend only on the permutation $\sigma(B)$, 
not on the product $B$ that effected $\sigma(B)$; in other words, both 
are independent of the switching history that effected the ultimate permutation.
These algorithms require the entries in $\sigma(B)$
to switch only with bodies in the set $\{x,y\}$, where  $x$ and $y$
are ``outsiders" who
did not participate in the initial mind-swapping spree.
Because these algorithms are designed for situations where
the switching history is forgotten,
it should come as no surprise that they may not be optimal in
situations where the switching history is known.
For example, 
Keeler's algorithm reverses $B_9$ 
with the product of thirteen transpositions given by
\begin{equation}\label{q13}
(x1)\cdot(y2)\cdot(x2)\cdot(y1)\cdot(x3)\cdot(x4)\cdot(x5)\cdot(x6)\cdot(x7)
\cdot(x8)\cdot(y9)\cdot(x9)\cdot(y3),
\end{equation}
while by \eqref{q9}, $B_9$ can be undone
with only nine transpositions.  We consider further comparisons of this kind
in Section 5.

We next look at the problem from {\em Stargate SG-1}.
For brevity, denote the bodies of Teal'c, Jack, Machello, and Daniel
by $1,2,3$, and $4$, respectively.  Recall that 
Daniel initially swapped minds with 
Machello, after which Teal'c swapped minds with Jack. This sequence of
two mind swaps can be represented by the product
\begin{equation}\label{h4}
H_4:=(12)\cdot(34).
\end{equation}
Samantha Carter observed that
$\sigma(H_4)=\sigma(Q_4)$, where
\begin{equation} \label{q4}
Q_4:=(24)\cdot(13)\cdot(23)\cdot(14).
\end{equation}
The products $H_4$ and $Q_4$ have no factor in common.
Thus $Q_4^{-1}$ serves to undo $H_4$ using four transpositions.
It follows easily from parity considerations (as in Lemma 3.3 below) that
$H_4$ cannot be undone with fewer than four transpositions, so that
$Q_4$ is best possible.

The product in \eqref{h4} is a special case of the product
\begin{equation} \label{h2r}
H_{2r}: = (12)\cdot(34)\cdots (2r-1,2r), \quad r \ge 1.
\end{equation}
Observe that
\begin{equation} 
\sigma(H_{2r}): = (12)(34)\cdots (2r-1,2r), \quad r \ge 1,
\end{equation}
since the factors of $H_{2r}$ are already disjoint.
We proceed to extend the problem of undoing $H_4$ to that of undoing $H_{2r}$.
The switching history now is a sequence of $r$ disjoint mind swaps
represented by \eqref{h2r}.
We are interested in finding an optimal product $Q$
for which $Q^{-1}$ undoes $H_{2r}$.  It turns out that such a product 
$Q$ needs no outsiders, except in the case $r=1$.

First consider the case $r=1$. Using the outsiders $3$ and $4$, we have
\begin{equation}
\sigma(Q_2)=\sigma(H_2)=(12), 
\end{equation}
where
\begin{equation}\label{q2}
Q_2:=(34)\cdot(23)\cdot(14)\cdot(13)\cdot(24).
\end{equation}
Since $(12)$ is not a factor of the product $Q_2$, it follows that
$Q_2^{-1}$ serves to undo $H_2$ using five transpositions.
We can show directly that
$H_2$ cannot be undone with fewer than five transpositions.
For suppose that $(12)$ were equal to
a product of $u<5$ distinct transpositions, none
equal to $(12)$.  Then since $(12)$
is an odd permutation,
we would have $u=3$ (again by Lemma 3.3), 
so that some product of four
distinct transpositions would equal the identity.  This is easily
seen to be impossible. 

Having already discussed the case $r=2$ in the paragraph containing 
\eqref{q4}, we next consider the case $r=3$.
Observe that
\begin{equation}
\sigma(Q_6)=\sigma(H_6)=(12)(34)(56), 
\end{equation}
where
\begin{equation} \label{q6}
Q_6:= (15)\cdot(25)\cdot(35)\cdot(46)\cdot(45)\cdot(16)\cdot(13).
\end{equation}
The product $Q_6$ has no factor in common with 
$H_6=(12)\cdot(34)\cdot(56)$.
Thus $Q_6^{-1}$ serves to undo $H_6$ using seven transpositions.

For general $r>3$,
we can undo $H_{2r}$ as follows.  If $r$ is even, undo the first two
factors of $H_{2r}$ using four transpositions as in \eqref{q4}, then repeat
the process for the next two factors of $H_{2r}$, etc. In this way, 
we undo $H_{2r}$ with $2r$ transpositions.  If $r$ is odd, 
undo the first three factors of $H_{2r}$
using seven transpositions as in \eqref{q6}, then undo consecutive
pairs of remaining factors as in the case where $r$ is even.  In this way, we
undo $H_{2r}$ with $2r+1$ transpositions.  In summary, for any $r \ge 2$, the
product $H_{2r}$ can be undone with $2r + \epsilon_r$ transpositions, where
\begin{equation}\label{epsilonr}
\epsilon_r =
\begin{cases}
\ 0, & \mbox{ if } \ 2\mid r  \
\\
\ 1, & \mbox{ if } \ 2\nmid r.
\end{cases}
\end{equation}
Is this best possible?  
Theorem 3.5 gives an affirmative answer for each $r \ge 2$.
In the exceptional case $r=1$, we've already seen
that $H_2$ can be undone with five transpositions, and this is best possible.

The proofs of Theorems 3.4 and 3.5 depend on Theorem 3.1, 
which may be of independent interest because it
answers the following question
about a permutation $\sigma$ expressed 
as a product of nontrivial disjoint cycles.
``When writing $\sigma$ as a product of distinct transpositions
none equal to a factor of $\sigma$, 
what is the smallest number of transpositions
that can be used?"

\section{Optimality theorems}

Let $\sigma$ be a permutation expressed as a product of $m \ge 1$ nontrivial 
disjoint cycles, and let $n\ge 2$ denote the number of entries in $\sigma$.
A paper by Mackiw \cite{mackiw} opens with this question:
When writing $\sigma$ as a product of transpositions, 
what is the smallest number of transpositions that can be used?
The well-known answer is $n-m$; see \cite{mackiw}, \cite[Theorem 4]{higgs}.

The following theorem answers a refinement of this question when $n >2$.
For brevity, if a transposition $\tau$ appears as a factor in the 
disjoint cycle factorization of $\sigma$, we will simply say
``$\tau$ occurs in $\sigma$".

\begin{thm}\label{thm3.1}
Let $\sigma$ be a permutation expressed as a product of $m \ge 1$ nontrivial
disjoint cycles, and let $n > 2$ denote the number of entries in $\sigma$.
Suppose that $\sigma=\sigma(Q)$, where $Q$ is a product of 
distinct transpositions, none occurring in $\sigma$.  Then
the smallest number $M$ of transposition factors that $Q$ can have is
\begin{equation}\label{M}
M = n-m+R+\epsilon_R,
\end{equation}
where $R$ is the number of transpositions occurring in $\sigma$,
and $\epsilon_R$ is as defined in \eqref{epsilonr}.
\end{thm}

It is interesting to note that $M$ depends only on
the number of transposition factors, the number of cycles, and the number
of entries in $\sigma$, but not on any additional information about the
cycle structure of $\sigma$.   For example, either 
$\sigma=(12)(345)(678)(9abcde)$ or
$\sigma=(1234)(56)(789a)(bcde)$ would yield the same value of $M$.

Theorem 3.1 does not hold for $n=2$, since then $M=5$ by
the argument following  \eqref{q2}.  It can be shown that Theorem 3.1
would be valid for all $n \ge 2$ if we dropped the requirement that the
factors of $Q$ be distinct.  In particular, for the case $n=2$, note that
$(12)=(23)(13)(23)$.

We prove Theorem 3.1 in the next section, using the following 
two well-known results from graph theory and group theory.
Proofs of these results may be found in
\cite[Theorem 11.2.1, p. 163]{cameron} and
\cite[pp. 82, 149]{BB}, respectively.

\begin{lem}\label{lem3.2}
A connected graph on $N$ vertices has at least $N-1$ edges.
\end{lem}

\begin{lem}\label{lem3.3}
If the product of $k$ transpositions in the symmetric group $S_n$
equals the identity, then $k$ must be even.
\end{lem}

Theorems 3.4 and 3.5 below
solve the optimality problems arising from
{\em Futurama} and {\em Stargate SG-1} that were posed in Section 2.
To prove Theorem 3.4, we will apply 
Theorem 3.1 with $\sigma$ chosen to be $\sigma(B_n)$. 
In this case, the lone
transposition $(12)$ occurring in $\sigma$ is also a factor of $B_n$.
To prove Theorem 3.5, we will apply
Theorem 3.1 with $\sigma$ chosen to be $\sigma(H_{2r})$.
In this case, all
transpositions occurring in $\sigma$ are factors of $H_{2r}$.
We are thereby able to exploit
the switching histories that effected $\sigma(B_n)$ 
and $\sigma(H_{2r})$,
even though no switching history is involved in Theorem 3.1 per se.

\begin{thm}\label{thm3.4}
Let $n \ge 5$.  Then
the product $B_n$ defined above \eqref{qn} cannot be undone with fewer than
$n$ transpositions.
\end{thm}

\begin{proof}
Recall from \eqref{dpn} that
$\sigma(B_n) = (12)(345\dots n)$. 
Suppose that
$\sigma(B_n) = \sigma(Q)$ for some product $Q$ of distinct 
transpositions, none equal
to a factor of $B_n$.  We need to show that $Q$ has at least $n$ factors.
By definition, $B_n$ has the factor $(12)$; thus
the product $Q$ cannot have $(12)$ as a factor.  It follows that
$Q$ is a product of distinct transpositions none
occurring in $\sigma(B_n)$.
We can therefore apply Theorem 3.1 with
$\sigma = \sigma(B_n)$ to conclude that the number of factors in $Q$ is at least
\[
M=n-m+R +\epsilon_R = n-2+1+1= n.  
\]
\end{proof}

\begin{thm}\label{thm3.5}
Let $r \ge 2$.  Then
the product $H_{2r}$ defined in \eqref{h2r} cannot be undone with fewer than
$2r+\epsilon_r$ transpositions.
\end{thm}

\begin{proof}
Suppose that
$\sigma(H_{2r}) = \sigma(Q)$ for some product $Q$ of distinct transpositions, 
none equal to a factor of $H_{2r}$.  
Therefore 
$Q$ is a product of distinct transpositions none
occurring in $\sigma(H_{2r})$.
Applying Theorem 3.1 with $\sigma=\sigma(H_{2r})$ and $R=r$,  
we conclude that the number
of factors in $Q$ is at least 
\[
M=n-m+r +\epsilon_r = 2r-r+r+\epsilon_r = 2r+\epsilon_r.  
\]
\end{proof}

\section{Proof of Theorem 3.1}

With $R=r$, we may express the permutation $\sigma$ 
in Theorem 3.1 as the following product
of disjoint cycles:
\begin{equation}\label{TH2r}
\sigma=\sigma(H_{2r})C_1\cdots C_{m-r},
\end{equation}
where $H_{2r}$ is defined in \eqref{h2r} and
the $C_i$ are disjoint cycles of length $\ell_i >2$.
(If $m=r$, then $\sigma$ is interpreted as $\sigma(H_{2r})$.)
Our object is to prove \eqref{M}.
We may suppose that $r \ge 1$, since if no transposition occurs in $\sigma$,
then $M=n-m$ as in \cite{mackiw}.

We proceed to show that $\sigma$ has parity $n-m +r +\epsilon_r$, 
by which we mean that the permutation
$\sigma$ is even or odd according as the integer $n-m +r +\epsilon_r$ 
is even or odd.
Each $C_i$ equals a product of $\ell_i-1$
transpositions; for example, $(abcdef) = (ab)(bc)(cd)(de)(ef)$.
Therefore $C_1\cdots C_{m-r}$ has parity
\[
\sum_{i=1}^{m-r} (\ell_i - 1) \ = \ (n-2r)-(m-r) = n - m -r.
\]
Since $\sigma(H_{2r})$ is a product of $r$ transpositions, 
it follows from \eqref{TH2r} that $\sigma$ has parity $n-m$.
Equivalently, $\sigma$ has parity $n-m +r + \epsilon_r$, since
$r + \epsilon_r$ is even by \eqref{epsilonr}.

Our object is to prove that
$M=n-m+r+\epsilon_r$ for $n>2$, $r \ge 1$.  We begin by showing that
\begin{equation} \label{upperbound}
M \le n-m+r+\epsilon_r.
\end{equation}
To accomplish this, it must be shown 
that $\sigma =\sigma(Q)$ for some product $Q$ of
$n-m+r+\epsilon_r$ distinct transpositions, none occurring in $\sigma$.

First consider the case $r=1$.
By \eqref{dpn}$-$\eqref{qn}, the permutation $(12)C_1$ equals
a product of $2+\ell_1$ distinct
transpositions, none equal to the transposition $(12)$.  
Thus $\sigma$ $=$ $(12)C_1\cdots C_{m-1}$
equals a product of
\[
2+\ell_1 + \sum_{i=2}^{m-1} (\ell_i -1) = n-m+2
\]
distinct transpositions, none equal to $(12)$.
This proves \eqref{upperbound} in the case $r=1$.  

Next let $r \ge 2$.  The statement containing \eqref{epsilonr}
tells us that for $r \ge 2$,
$\sigma(H_{2r})$ equals a product of
$2r+\epsilon_r$ distinct transpositions, none equal to a factor of $H_{2r}$.
Thus  $\sigma=\sigma(Q)$ for a product $Q$ of
\[
2r+\epsilon_r + \sum_{i=1}^{m-r} (\ell_i -1) = n-m+r+\epsilon_r
\]
distinct transpositions, none occurring in $\sigma$.
This completes the proof of \eqref{upperbound}.
(Note that $Q$ has been constructed such that
every entry in $Q$ is an entry in $\sigma$.
This is in contrast to the case $n=2$, as we see from \eqref{q2}).

It remains to prove the reverse inequality 
\begin{equation} \label{lowerbound1}
M \ge  n-m+r+\epsilon_r.
\end{equation}
We need only prove
\begin{equation} \label{lowerbound}
M \ge  n-m+r;
\end{equation}
indeed, \eqref{lowerbound1} and \eqref{lowerbound} are equivalent 
by Lemma 3.3, since as was noted above, $\sigma$ 
has parity $n-m +r  +\epsilon_r$.
Supposing that $\sigma=\sigma(Q)$ for a product $Q$ of $w$
distinct transpositions, none occurring in $\sigma$, our goal
is now to show that $w \ge n-m+r$.  

To illuminate the exposition, we
will use the running example
\[
\sigma=(12)(34)(56)(78)(9ab)(cde) = \sigma(Q), 
\]
with $Q$ equal to the product
\[
(15)\cdot(89)\cdot(25)\cdot(7b)\cdot(4x)\cdot(35)\cdot(7a)\cdot(46)\cdot(79)
\cdot(45)\cdot(8b)\cdot(16)\cdot(13)\cdot(cd)\cdot(de)\cdot(3x),
\]
where $a,b,c,d,e,x$ stand for $10,11,12,13,14,15$.
In this example, $n=14$, $r=4$, $m=6$, and $w=16$.

Let $G$ be a graph whose vertex set $V(G)$ is the set of entries in $Q$,
and whose edges $[ij]$ correspond to the $w$ factors $(ij)$ of $Q$.
The graph $G$ for the product $Q$ in our example is illustrated in Figure 1.
Since $\sigma=\sigma(Q)$,
the set $V(G)$ contains all the entries in $\sigma$, but
$V(G)$ may also contain entries in $Q$ outside of $\sigma$.
In our example, $x$ is the only entry in $Q$ outside of $\sigma$.

For each $i$ with $1\le i \le m-r$, $G$ has a connected component
$G_i$ such that $V(G_i)$ contains all of the entries in the cycle $C_i$.
Let $J$ denote the union of the components $G_i$.  (View $J$ as empty
if $m=r$.)  In our example, $J$ is the union of two components
$G_1$, $G_2$ with $V(G_1) = \{7,8,9,a,b\}$, $V(G_2) = \{c,d,e\}$,
so $V(J)= \{7,8,9,a,b,c,d,e\}$.

The set
$\{1,2,\dots, 2r\} \subset V(G)$ can be written as a disjoint union
\begin{equation}
\{1,2,\dots, 2r\} = A \cup Z,
\end{equation}
where $A \subset V(J)$ and $Z$ is disjoint from $V(J)$. 
For the cardinalities, write
$\alpha:=|A|$, so that $|Z| =2r -\alpha$.
In our example, $A=\{7,8\}$ with $\alpha=2$ and
$Z=\{1,2,3,4,5,6\}$.

Since $V(J)$ contains the $n-2r$ entries in the cycles $C_i$, we have
$|V(J)| \ge n - 2r + \alpha$.  By Lemma 3.2,
each $G_i$ contains at least $|G_i|-1$ edges.  
Therefore, since the graph $J$ is the disjoint union of at most $m-r$ 
distinct components $G_i$, the number of edges in $J$ is at least
\[
(n - 2r + \alpha) - (m-r) = n-m +\alpha -r .
\]
It follows that $Q$ has at least $n-m + \alpha -r$ distinct transposition
factors whose entries are all in $V(J)$.  In particular,
$w \ge  n-m +\alpha -r$, but this does not quite accomplish our goal
of showing that $w \ge n-m+r$. Thus we seek additional transposition
factors of $Q$,
this time with entries {\em outside} of $V(J)$.

With the aim of obtaining these additional factors of $Q$, 
we will construct a set $E$ with $Z \subset E \subset V(G)$
such that $E$ is disjoint from $V(J)$ and
the elements of $E$ fill at least $2|E|$ slots in $Q$.
(Note that any element of $E$ not in $Z$ is necessarily 
outside of $\sigma$.)  Since $|E| \ge |Z|$, we will  obtain
the desired lower bound
\[
w \ge (n-m + \alpha -r)  + |E|  \ge   (n-m + \alpha -r)  + |Z|  = n-m+r.
\]

If each of the elements in $Z$ fills at least two slots in $Q$,
then we may take $E=Z$ and the proof is complete.
It remains to construct $E$ in the case 
where $Z$ contains ``singletons," i.e.,
elements that occur only once as entries in $Q$.
Note that a singleton corresponds to a vertex of degree 1 in $V(G)$.
In our example, $2$ is the only singleton in $Z=\{1,2,3,4,5,6\}$.

We proceed to specify the elements that will make up the desired set
$E$ defined in \eqref{E} below.
A singleton $b_1 \in Z$ is paired in $Q$ with some $g \in V(G)$, that is,  
$Q$ has a factor $(gb_1)$. In our example,  $(gb_1)=(25)$.
Note that $g$ cannot be a singleton, otherwise $(gb_1)$ would be a factor
of $Q$ that occurs in $\sigma$.
Since $b_1 \notin V(J)$, we have $g \notin V(J)$.  
Let $S(g):=\{b_1,\dots, b_k \}$ be the
set of all singletons in $Z$ that are paired with $g$, and let
\[
(gb_1), \ (gb_2), \ \dots , \ (gb_k)
\]
be transpositions appearing in that same left-to-right order 
among the factors of $Q$.
In our example, $k=1$.

Since no factor of $Q$ occurs in $\sigma$ and $b_1 \in Z$, 
the permutation $\sigma$
cannot map the singleton $b_1$ to $g$.  Thus $g$ occurs as an entry in
$Q$ to the left of the factor $(gb_1)$.  Similarly, $\sigma$ cannot map $g$
to the singleton $b_k$, so $g$ occurs as an entry in $Q$ to the right of
$(gb_k)$.  Therefore $g$ occurs at least $k+2$ times as an entry,
so that all together, the $k+1$ elements in 
\[
U(g):=\{b_1, \dots, b_k, g\} = S(g) \cup \{g\}
\]
fill at least $2(k+1)=2|U(g)|$ slots in $Q$.  In our example, $(15)$ occurs
to the left of the factor $(gb_1)=(25)$ in $Q$, both $(35)$ and $(45)$
occur to the right of $(25)$, and the two elements in $U(5)=\{2,5\}$
together fill five slots in $Q$. 

If $Z$ contains another singleton $b_1' \notin S(g)$, 
then repeat this procedure with an element $g' \in V(G)$ 
that is paired with $b_1'$.  Clearly $g' \neq g$.
Note that $S(g')$ is disjoint from $S(g)$, since an element common to both
sets would have to occur at least twice as an entry in $Q$.
Thus $U(g')$ is disjoint from $U(g)$.
Repeat the procedure again and again until all the singletons in $Z$ have 
been exhausted.  

Let $Z^*$ denote the set of all non-singletons in $Z$ which are not
in the set $\{g, g', \dots\}$.  The  elements of
$Z^*$ fill at least $2|Z^*|$ slots in $Q$.  Define $E$ to be the
disjoint union
\begin{equation}\label{E}
E = Z^* \cup  U(g) \cup U(g') \cup \cdots .
\end{equation}
In our example, $E=Z$ (but $E$ would strictly contain $Z$ in an example where
the set $\{g, g', \dots\}$ contained an element outside of $\sigma$).  
The set $E$ in \eqref{E} satisfies the required conditions, since
\[
E = Z \cup \{g, g', \dots\}\subset V(G), 
\]
$E$ is disjoint from $V(J)$, and the elements of $E$
fill at least $2|E|$ slots in $Q$. \qed

\section{Conclusion and related problems}

Let $\sigma(B)$ be a known permutation in $S_n$
corresponding to a product
$B$ representing a (possibly forgotten) switching history.
This section begins with a summary of optimal ways 
for undoing $B$ when 
\begin{enumerate}[(i)]
\item $B$ is unknown,
\item  $B=H_{2r}$ (defined in \eqref{h2r}),
\item  $B=B_n$ (defined above \eqref{qn}).
\end{enumerate}
We then discuss the situation for some other classes of products $B$.
Our examples illustrate how dramatically different 
the optimal solutions can be
for products $B$ and $B'$ representing different
switching sequences, even when $B$ and $B'$ effect the same permutation
$\sigma(B) = \sigma(B')$.

Assume that $\sigma(B)$ has $n$ entries and
$m$ nontrivial disjoint cycles.  Let $j(B)$ denote the number of
factors in an optimal product $Q$ for which $Q^{-1}$ undoes $B$.
If $B$ represents a {\em forgotten} switching history, 
then by \cite[Theorem 1]{keeler},
\begin{equation}\label{bp}
j(B) = n+m+2,
\end{equation}
and $Q$ requires two outsiders.
In cases where there is full information about $B$, we can often reduce
the right side of \eqref{bp}.  
For example, if $B=H_n=H_{2r}$ with $r \ge 2$, then
Theorem 3.5 and the statement containing \eqref{epsilonr} show that
\[
j(B)=n\ \ \mbox{or} \ \  j(B)=n+1 \ \ \mbox{according as}  
\ r \  \mbox{is even or odd},
\]
and no outsiders are needed.
Even if only partial information is available about $B$, for instance when
$B=B_n$, we may be able to reduce the right side of \eqref{bp}.
For example,  Theorem 3.4 and \eqref{qn} show that
\begin{equation}\label{bnp}
j(B_n)=n,  \quad \mbox{for} \ \ n \ge 5,
\end{equation}
and no outsiders are needed.

Let $B_n'$ denote the product $Q_n$ defined in \eqref{qn}.
Then $\sigma(B_n')= \sigma(B_n)$.
We can reduce the right side of
\eqref{bnp} upon replacing $B_n$ by $B_n'$, as follows:
\begin{equation}\label{bnp1}
j(B_n')=n-2,  \quad \mbox{for} \ \ n \ge 5.
\end{equation}
The idea behind \eqref{bnp1} is 
to undo $B_n'$ with the product 
\[
(12)\cdot(3n)\cdot(4n)\cdots (n-1,n).
\]
One can show that this way of undoing $B_n'$ is best possible, 
by using Lemma 3.3 and
\cite[Lemma 1(i),(ii)]{keeler}.
For an example with $n=5$, we have
\[
B_5' = (23)\cdot(15)\cdot(14)\cdot(13)\cdot(25) \ \ \mbox{and} \ \ 
\sigma(B_5') = \sigma(B_5)=(12)(345),
\]
and $\ j(B_5') = 3$, since $(12)(35)(45)$ undoes $B_5'$.

Yet another product $B$ with $\sigma(B)=(12)(345)$
is 
\[
B=(25)\cdot(45)\cdot(35)\cdot(14)\cdot(24). 
\]
Since $\sigma(B)=\sigma(Q)$ for $Q:=(23)\cdot(12)\cdot(15)\cdot(13)\cdot(34)$,
$Q^{-1}$ serves to undo $B$ with five transpositions. 
We leave it as an exercise for the reader to show that 
this is best possible, i.e., $j(B)=5$.

We proceed to give the values of $j(B)$ for three
different classes of products $B$ which effect
the same $n$-cycle $\ \sigma(B)=(123\dots n)$.
Again these examples illustrate how sensitive the optimal solutions
can be to the switching history.
It would be interesting to find $j(B)$ for other natural
classes of products $B$, and to determine how many outsiders are
required for undoing these $B$.

(a) If 
\[
B=(12)\cdot(23)\cdot(34)\cdots (n-1,n) \ \ \mbox{with} \ \  n \ge 5, 
\]
then 
$\sigma(B)=(123\dots n)$ and $B$ can be undone with $n+1$
transpositions, without outsiders.  Moreover, 
by \cite[Theorem 2]{keeler},
$n+1$ is best possible, i.e., $j(B)=n+1$.

(b) If 
\[
B=(n,n-1)\cdots (n3)\cdot(n2)\cdot(n1) \ \ \mbox{with} \ \  n \ge 3, 
\]
then
$\sigma(B)=(123\dots n)$ and $B$ can be undone with $n+1$
transpositions, where one outsider is necessary and sufficient.  Moreover,
by \cite[Theorem 3]{keeler},
$n+1$ is best possible, i.e., $j(B)=n+1$.

(c) In (a) and (b) above, $B$ has only $n-1$ factors.
If $B$ is a product of all $(n^2-n)/2$  transpositions in $S_n$,
then two outsiders are necessary and sufficient to undo $B$, and
by \cite[Theorem 1]{keeler},
$j(B) = n+m+2$.  As shown in \cite{keeler}, when say
$n$ is congruent to 1 modulo 4, it is 
not hard to give a recipe for constructing
such $B$ in $S_n$ with $\ \sigma(B)=(123\dots n)$.
For an example in $S_5$,
consider the product of all ten transpositions
\[
B=(54)\cdot(53)\cdot(52)\cdot(51)\cdot(12)\cdot(23)\cdot(14)\cdot(13)
\cdot(24)\cdot(34);
\]
then $\sigma(B)=(12345)$ and $j(B)=8$.

Examples of two products $B$ and $B'$ which effect the same
4-cycle $\sigma(B)=\sigma(B') = (1234)$
are
\[
B=(24)\cdot(23)\cdot(14) \ \ \mbox{and} \ \  
B'=(23)\cdot(12)\cdot(34)\cdot(14)\cdot(24).
\]
It is easy to see that $(34)\cdot(12)\cdot(13)$ undoes $B$ and $j(B)=3$.
We close with another fun problem for the reader:
show that two outsiders are necessary and sufficient to undo $B'$,
and $j(B') = 7$.

\section*{Acknowledgements}
The authors are very grateful to the referees and to the editor for many 
helpful suggestions.

\noindent \textbf{Ron Evans} received a BS from the University of Michigan 
in 1967 and a Ph.D. from the University of Illinois in 1974. After a two-year 
instructorship at the University of Wisconsin, he joined the UCSD
Mathematics Department in 1975. His main interest is in number theory, 
with recent work focusing on character sum evaluations with connections 
to Hecke eigenfunctions.  He is coauthor of {\em Gauss and Jacobi sums} 
(Wiley, 1998) with his Ph.D. advisor Bruce C. Berndt and Kenneth S. Williams.

\vspace{5 mm}

\noindent \textbf{Lihua Huang} obtained her BS in mathematics from the
University of California at San Diego in 2012.
Her research interests include group theory and several complex variables.
She has presented her research at conferences 
held by Ohio State, Cal Poly Pomona, Penn State, UCSD, and the 2013 
Joint Mathematics Meetings in San Diego.  In 2012, she became a graduate 
student in mathematics at Rutgers University, where she received a GAANN 
Fellowship. 
\end{document}